\documentclass[12pt,draft]{amsart}
\usepackage{amsmath,amsthm,latexsym,amscd,amsbsy,amssymb,pb-diagram}
\chardef\bslash=`\\ 

\makeatletter
\def\verbatim{\interlinepenalty\@M \@verbatim
  \leftskip\@totalleftmargin\advance\leftskip2pc
  \frenchspacing\@vobeyspaces \@xverbatim}
\makeatother
\makeatletter
  \def\dgt@k{\dg@DX=-3 \dg@DY=2 \dg@SIZE=3} 
\makeatother

\makeatletter
  \def\dgt@kk{\dg@DX=3 \dg@DY=-1 \dg@SIZE=3}%
\makeatother

\theoremstyle{plain}
\newtheorem{thm}{Theorem}[section]
\newtheorem{cor}[thm]{Corollary}

\newtheorem{pro}[thm]{Proposition}

\theoremstyle{definition}

\newtheorem{defin}[thm]{Definition}

\numberwithin{equation}{section}

\newcounter{rmnum}

\def\symbolnote#1#2{\let\thefootn=\thefootnote%
\renewcommand{\thefootnote}{\fnsymbol{footnote}}%
\footnotemark[#1]%
\footnotetext[#1]{#2}%
\let\thefootnote=\thefootn
}

\newfont{\bbb}{msbm10 scaled \magstep1}
\newfont{\bbc}{msbm8 scaled \magstep0}

\newcommand{\N}{\mbox{\bbb N}}

\begin{document}


\title[$C^{\ast}$-algebras of infinite real rank]
{$C^{\ast}$-algebras of infinite real rank}
\author{A.~Chigogidze}
\address{Department of Mathematics and Statistics,
University of Saskatche\-wan,
McLe\-an Hall, 106 Wiggins Road, Saskatoon, SK, S7N 5E6,
Canada}
\email{chigogid@math.usask.ca}
\thanks{The first author was partially supported by NSERC research grant.}
\author{V.~Valov}
\address{Department of Mathematics and Computer Science, Nipissing University,
100 College Drive, P.O. Box 5002, North Bay, ON, P1B 8L7, Canada}
\email{veskov@unipissing.ca}
\thanks{The second author was partially supported by Nipissing University Research Council Grant.}


\begin{abstract}{We introduce the notion of weakly (strongly) infinite real rank for unital $C^{\ast}$-algebras. It is shown that a compact space $X$ is weakly (strongly) infine-dimensional if and only if $C(X)$ has weakly (strongly) infinite real rank. Some other properties of this concept are also investigated. In particular, we show that the group $C^{\ast}$-algebra $C^{\ast}\left( {\mathbb F}_{\infty}\right)$ of the free group on countable number of generators has strongly infinite real rank.}
\end{abstract}

\maketitle
\markboth{A.~Chigogidze, V.~Valov}{$C^{\ast}$-algebras of infinite real rank}


It is clear that some $C^{\ast}$-algebras of infinite real rank  have infinite rank in a very strong sense of this word, while others do not. In order to formally distinguish these types of infinite ranks from each other we introduce the concept of weakly (strongly) infinite real rank.
Proposition \ref{P:definitions} characterizes usual real rank in terms of {\em infinite} sequences of self-adjoint elements and serves as a basis of our definition \ref{D:weakly}. We completely settle the commutative case by proving (Theorem \ref{T:iff}) that the algebra $C(X)$ has weakly infinite real rank if and only if $X$ is a weakly infinite dimensional compactum. As expected, the group
$C^{\ast}$-algebra $C^{\ast}\left( {\mathbb F}_{\infty}\right)$ 
of the free group on countable number of generators
has strongly infinite real rank (Corollary \ref{C:notweakly}).

\section{Preliminaries}\label{S:pre}
All $C^{\ast}$-algebras below are assumed to be unital. The set of all
self-adjoint elements of a $C^{\ast}$-algebra $X$ is denoted
by $X_{sa}$. 

The real rank of a unital $C^{\ast}$-algebra $X$, denoted by $\operatorname{rr}(X)$,
is defined as follows \cite{brownped91}. We say that $\operatorname{rr}(X) \leq n$ if for each
$(n+1)$-tuple $(x_{1},\dots ,x_{n+1})$ of self-adjoint elements
in $X$ and every $\epsilon > 0$,
there exists an $(n+1)$-tuple $(y_{1},\dots ,y_{n+1})$ in $X_{sa}$
such that $\sum_{k=1}^{n+1} y_{k}^{2}$ is invertible and
$\left\| \sum_{k=1}^{n+1} (x_{k}-y_{k})^{2}\right\| < \epsilon$ .


\subsection{Alternative definitions of the real rank}\label{SS:adrr}
It is interesting that the real rank can be equivalently defined in terms of
infinite sequences.

\begin{pro}\label{P:definitions}
Let $X$ be a unital $C^{\ast}$-algebra. Then the following conditions are equivalent:
\begin{itemize}
\item[(i)]
$\operatorname{rr}(X) \leq n$.
\item[(ii)]
for each $(n+1)$-tuple $(x_{1},\dots ,x_{n+1})$ 
in $X_{sa}$ and for each $\epsilon > 0$,
there exists an $(n+1)$-tuple $(y_{1},\dots ,y_{n+1})$ in $X_{sa}$
such that $\sum_{k=1}^{n+1} y_{k}^{2}$ is invertible and
$\left\| x_{k}-y_{k}\right\| < \epsilon$
for each $k = 1,2,\dots ,n+1$.
\item[(iii)]
for any sequence of self-adjoint elements $\{x_i \colon i \in N\}\subseteq X_{sa}$
and for any sequence of positive real numbers $\{\epsilon_i\colon i\in\N\}$ there exists a sequence
$\{y_i \colon i \in \N\}\subseteq X_{sa}$ such that 
\begin{itemize}
\item[(a)]
$\|x_i-y_i\|<\epsilon_i$, for each $i\in\N$, 
\item[(b)]
for any subset $D \subseteq \mathbb{N}$, with $| D| = n+1$, the element
$\sum _{i\in D}y_{i}^2$ is invertible.
\end{itemize}
\item[(iv)]
for any sequence of self-adjoint elements $\{x_i \colon i \in \N\}\subseteq X_{sa}$
and for any $\epsilon >0$ there exists a sequence
$\{y_i \colon i \in \N\}\subseteq X_{sa}$ such that 
\begin{itemize}
\item[(a)]
$\|x_i-y_i\|<\epsilon$, for each $i\in\N$, 
\item[(b)]
for any subset $D \subseteq \mathbb{N}$, with $| D| = n+1$, the element
$\sum _{i\in D}y_{i}^2$ is invertible.
\end{itemize}
\item[(v)]
for any sequence of self-adjoint elements $\{x_i \colon i \in \N\}\subseteq X_{sa}$
such that $\| x_{i}\| = 1$ for each $i \in \mathbb{N}$ and
for any $\epsilon >0$ there exists a sequence
$\{y_i \colon i \in \N\}\subseteq X_{sa}$ such that 
\begin{itemize}
\item[(a)]
$\|x_i-y_i\|<\epsilon$, for each $i\in\N$, 
\item[(b)]
for any subset $D \subseteq \mathbb{N}$, with $| D| = n+1$, the element
$\sum _{i\in D}y_{i}^2$ is invertible.
\end{itemize}
\end{itemize} 
\end{pro}
\begin{proof}
(i) $\Longrightarrow$ (ii). Let  $(x_{1},\dots ,x_{n+1})$  be an
$(n+1)$-tuple in $X_{sa}$ and $\epsilon > 0$. By (i), there exists an $(n+1)$-tuple $(y_{1},\dots ,y_{n+1})$ in $X_{sa}$
such that $\sum_{k=1}^{n+1} y_{k}^{2}$ is invertible and
$\left\| \sum_{k=1}^{n+1} (x_{k}-y_{k})^{2}\right\| < \epsilon^{2}$.

Since $x_{k} -y_{k} \in X_{sa}$, it follows (\cite[2.2.4 Theorem]{murphy}) that $(x_{k}-y_{k})^{2} \geq 0$ for each $k = 1,\dots ,n+1$. Then, by \cite[2.2.3 Lemma]{murphy},
$\sum_{k=1}^{n+1} (x_{k}-y_{k})^{2} \geq 0$. Note also that 
$\sum_{i=1}^{n+1} (x_{i}-y_{i})^{2} - (x_{k} -y_{k})^{2} = \sum_{i=1,i \neq k}^{n+1} (x_{i}-y_{i})^{2} \geq 0$, $k = 1,\dots ,n+1$, which guarantees that $(x_{k} -y_{k})^{2} \leq \sum_{i=1}^{n+1} (x_{i}-y_{i})^{2}$ for each $k = 1,\dots ,n+1$. By \cite[2.2.5 Theorem]{murphy}, 
$\| x_{k}-y_{k}\|^{2}  = \| (x_{k}-y_{k})^{2}\| \leq \left\| \sum_{i=1}^{n+1} (x_{i}-y_{i})^{2}\right\| < \epsilon^{2}$. Consequently, $\| x_{k}-y_{k}\| < \epsilon$, $k = 1,\dots ,n+1$. This shows that condition (ii) is satisfied.

(ii) $\Longrightarrow$ (iii). Suppose
$\operatorname{rr}(X) \leq n$ and let $\{x_i \colon i \in \N\}\subset X_{sa}$ and $\{\epsilon_i \colon i \in \N\}$ 
be sequences of self-adjoint elements of $X$ and positive real
numbers, respectively. Denote by
$\mathbb{N}_{n+1}$ the family of all subsets of $\mathbb{N}$ of
cardinality $n+1$. For every $i\in\mathbb{N}$ and
$D\in \mathbb{N}_{n+1}$ let
$H_i=\{x\in X_{sa}:\|x-x_i\|\leq 2^{-1}\epsilon_i\}$ and
$H_D=\prod\{H_i:i\in D\}$ be the topological product of all $H_{i}$, $i \in D$. We also consider the topological product
$H=\prod\{H_i:i\in\N\}$ and the natural projections
$\pi_D\colon H\to H_D$. Define the continuous maps
$\phi_D\colon X_{sa}^D\to X$, 
$\phi_D(z_{i_1},z_{i_2},..,z_{i_{n+1}})=
\sum_{j=1}^{n+1}z_{i_j}^2$, $D\in \mathbb{N}_{n+1}$.
Since the real rank of $X$ is $n$, Lemma 3.1(ii) yields
that $\phi_D^{-1}(G)$ is dense (and, obviously, open)
in $X_{sa}^D$ for every $D\in \mathbb{N}_{n+1}$, where
$G$ is the set of all invertible elements of $X$. The last
observation implies that each set
$G_D=\phi_D^{-1}(G)\cap H_D$ is open and dense in $H_D$.
Consequently, each $U_D=\pi_D^{-1}(G_D)$ is open and dense
in $H$ because the projections $\pi_D$ are continuous and
open maps. Finally, using that $H$ (as a product of countably
many complete metric spaces) has the Baire property, we
conclude that the intersection $U$ of all $U_D$ is non-empty.
Take any point $(y_i)$ from $U$. Then $y_i\in H_i$, so
$y_i\in X_{sa}$ and $\|x_i-y_i\|\leq 2^{-1}\epsilon_i<\epsilon_i$
for every $i\in\mathbb{N}$. Moreover, for any
$D\in \mathbb{N}_{n+1}$ the point $y_D=(y_i)_{i\in D}$
belongs to $G_D$, hence  $\sum\{y_i^2:i\in D\}$ is invertible.   

Implications (iii) $\Longrightarrow$ (iiv) and (iv) $\Longrightarrow$ (v)
are trivial.

(v) $\Longrightarrow$ (i). Let $(x_{1},\dots ,x_{n+1})$ be an $(n+1)$-tuple
of non-zero self-adjoint elements in $X$ and $\epsilon > 0$. Consider the 
sequence $\{ \bar{x}_{i}\}$ of self-adjoint elements of $X$, where
\[ \bar{x}_{i} =
\begin{cases}
\displaystyle \frac{x_{i}}{\| x_{i}\|},\;\;\text{if}\;\; i \leq n+1 ;\\
1 ,\;\; \text{if}\;\; i > n+1 .
\end{cases}
\] 

By (iv), there exists a sequence $\{ \bar{y}_{i}\colon i \in \N\}$ of self-adjoint elements of $X$ such that
$\sum_{i=1}^{n+1}\bar{y}_{i}^{2}$ is invertible and $\displaystyle\| \bar{x}_{i} - \bar{y}_{i}\| < \frac{\epsilon}{\max\{\| x_{i}\| \colon i = 1,\dots n+1\}}$ for each $i \in \mathbb{N}$. Now let $y_{i} = \| x_{i}\|\cdot \bar{y}_{i}$, $i \in \mathbb{N}$.
Then for every $i \leq n+1$ we have
\begin{multline*}
 \| x_{i} - y_{i}\| = \left\Vert \| x_{i}\|\cdot \bar{x}_{i} - \| x_{i}\|\cdot \bar{y}_{i}\right\Vert = \| x_{i}\| \cdot \| \bar{x}_{i} - \bar{y}_{i}\| <\\ \| x_{i}\| \cdot \frac{\epsilon}{\max\{\| x_{i}\| \colon i = 1,\dots n+1\}} \leq \epsilon .
\end{multline*}

The invertibility of $\sum_{i=1}^{n+1}\bar{y}_{i}^{2}$ is equivalent to the validity of the equation $1 = \sum_{i=1}^{n+1}z_{i}\bar{y}_{i}$ for a suitable $(n+1)$-tuple $(z_{1},\dots ,z_{n+1})$. Clearly $\displaystyle 1 = \sum_{i=1}^{n+1}\frac{z_{i}}{\| x_{i}\|}\cdot \| x_{i}\| \bar{y}_{i} = \sum_{i=1}^{n+1}\frac{z_{i}}{\| x_{i}\|}\cdot  y_{i}$ which in turn implies the invertibility of $\sum_{i=1}^{n+1}y_{i}^{2}$. 
\end{proof}


\subsection{Bounded rank}\label{SS:br}
For the readers convenience below we present definitions and couple of results related to the bounded rank. Details can be found in \cite{chivabr}.

\begin{defin}\label{D:modifiedunessential}
Let $K >0$. We say that an $m$-tuple $( y_{1},\dots ,y_{m})$ of self-adjoint elements of a
unital $C^{\ast}$-algebra $X$ is $K$-{\em unessential} if for every rational $\delta > 0$ there exists
an $m$-tuple $(z_{1},\dots ,z_{m})$ of self-adjoint elements of $X$ satisfying the following conditions:
\begin{itemize}
\item[(a)]
$\| y_{k}-z_{k}\| \leq\delta$ for each $k = 1,\dots ,m$,
\item[(b)]
The element $\sum_{k=1}^{m}z_{k}^{2}$ is invertible and $\left\| \left(\sum_{k=1}^{m}z_{k}^{2}\right)^{-1}\right\| \leq \displaystyle\frac{1}{K\cdot \delta^{2}}$.
\end{itemize}
$1$-unessential tuples are referred as {\em unessential}.
\end{defin}

\begin{defin}\label{D:modifiedbrr}
Let $K>0$. We say that the {\em bounded rank} of a unital $C^{\ast}$-algebra $X$ with respect to $K$ 
does not exceed $n$ (notation: $\operatorname{br}_{K}(X) \leq n$) if for any 
$(n+1)$-tuple $(x_{1},\dots ,x_{n+1})$ of self-adjoint elements of $X$ and for any $\epsilon > 0$ there exists a $K$-unessential $(n+1)$-tuple $(y_{1},\dots ,y_{n+1})$ in $X$ such that $\| x_{k}-y_{k}\| < \epsilon$ for each $k = 1,\dots ,n+1$. For simplicity $\operatorname{br}_{1}(X)$  is denoted by $\operatorname{br}(X)$ and it is called a {\em bounded rank}.
\end{defin}

\begin{pro}\label{P:commutting}
Let $(y_{1},\dots ,y_{m})$ be a commuting $m$-tuple of self-adjoint elements of the unital $C^{\ast}$-algebra $X$. If $\sum_{i=1}^{m}y_{i}^{2}$ is invertible, then $(y_{1},\dots ,y_{m})$ is $K$-unessential for any positive $K\leq 1$.
\end{pro}

\begin{cor}\label{C:rrcommutative}
Let $X$ be a commutative unital $C^{\ast}$-algebra and $0<K\leq 1$. Then $\operatorname{br}_{K}(X) = \operatorname{rr}(X) = \dim \Omega (X)$, where $\Omega (X)$ is the spectrum of $X$.
\end{cor}


\section{Infinite rank}\label{S:weakly}

We begin by presenting the definition of weakly infinite real rank.

\subsection{Weakly (strongly) infinite real and bounded ranks}\label{SS:weakly}

\begin{defin}\label{D:weakly}
We say that
a unital $C^*$ algebra $X$ has a weakly infinite real rank if for
any sequence of self-adjoint elements $\{x_i \colon i \in \mathbb{N}\}\subset X_{sa}$
and any $\epsilon > 0$ there is a sequence
$\{y_i \colon i \in \mathbb{N}\}\subset X_{sa}$ such that $||x_i-y_i||<\epsilon$ for every  
$i\in\mathbb{N}$ and the element
$\sum_{i\in D}y_i^2$ is invertible for some finite set $D$ of indices. If $X$ does not have weakly infinite real rank, then we say that $X$ has strongly infinite real rank.
\end{defin}

The bounded version can be defined similarly. 

\begin{defin}\label{D:sequnessential}
Let $K > 0$. A sequence of self-adjoint elements of a unital $C^{\ast}$-algebra is $K$-{\em unessential} if it contains a finite $K$-unessential (in the sense of Definition \ref{D:modifiedunessential}) subset.
\end{defin}

\begin{defin}\label{D:bweakly}
Let $K > 0$. We say that
a unital $C^*$ algebra $X$ has a weakly infinite bounded rank with respect to $K$ if for
any sequence of self-adjoint elements $\{x_i \colon i \in \mathbb{N}\}\subset X_{sa}$
and any $\epsilon > 0$ there is a $K$-unessential sequence
$\{y_i \colon i \in \mathbb{N}\}\subset X_{sa}$ such that $||x_i-y_i||<\epsilon$ for every  
$i\in\mathbb{N}$. If $X$ does not have weakly infinite bounded rank, then we say that $X$ has strongly infinite bounded rank.
\end{defin}

For future references we record the following statement.

\begin{pro}\label{P:finite}
Every unital $C^{\ast}$-algebra of a finite real rank has weakly infinite real rank.
\end{pro}
\begin{proof}
Apply Proposition \ref{P:definitions}.
\end{proof}

Note that, as it follows from Proposition \ref{P:example}, the converse of Proposition \ref{P:finite} is not true.

\begin{pro}\label{L:homomorphism}
Let $f \colon X \to Y$ be a surjective $\ast$-homomorphism of
unital $C^{\ast}$-algebras. If $X$ has  weakly infinite real rank, then so does $Y$.
\end{pro}
\begin{proof}
For any sequence of self-adjoint elements $\{y_i \colon i \in \N \}\subset Y_{sa}$
and for any $\epsilon > 0$ we need to find a sequence
$\{z_i \colon i \in \N\}\subset Y_{sa}$ such that 
\begin{itemize}
\item[(i)$_{Y}$]
$||y_i-z_i||<\epsilon$ for every $i\in\mathbb{N}$,
\item[(ii)$_{Y}$]
for some $k\geq 1$ the element
$\sum _{i=1}^{k}z_i^2$ is invertible. 
\end{itemize}
For every $i\in\mathbb{N}$ let $x_{i}$ be a self-adjoint element
in $X$ such that $f(x_{i}) = y_{i}$. Since $X$ has weakly infinite
real rank, there exists a
sequence $\{ w_{i} \colon i \in\N \}$ of self-adjoint elements in $X$ such that
\begin{itemize}
\item[(i)$_{X}$]
$\| x_i-w_i\| < \epsilon$ for every $i\in\mathbb{N}$,
\item[(ii)$_{X}$]
the element
$\sum _{i=1}^{k}w_i^2$ is invertible for some $k \geq 1$. 
\end{itemize}

Let $z_{i} =f(w_{i})$, $i \in \mathbb{N}$. Clearly $z_{i} \in Y_{sa}$ and
$\| y_{i}-z_{i}\| = \| f(x_{i}) - f(w_{i})\| = \| f(x_{i}-w_{i})\| \leq \| x_{i}-w_{i}\| < \epsilon$ for each $i \in \mathbb{N}$. By (ii)$_{X}$, there exists an element $a \in X$ such that $a\cdot \sum _{i=1}^{k}w_i^2 = 1$. Clearly $f(a)\cdot \sum _{i=1}^{k}z_i^2 = f(a)\cdot \sum _{i=1}^{k}f(w_i)^2 = f(a)\cdot f\left(\sum _{i=1}^{k}w_i^2 \right)= f\left( a\cdot \sum _{i=1}^{k}w_i^2\right) = f(1) = 1$ which shows that $\sum _{i=1}^{k}z_i^2$ is invertible.
\end{proof}

\begin{pro}\label{L:bhomomorphism}
Let $K > 0$ and $f \colon X \to Y$ be a surjective $\ast$-homomorphism of
unital $C^{\ast}$-algebras. If $X$ has  weakly infinite bounded rank with respect to $K$, then so does $Y$.
\end{pro}
\begin{proof}
For a sequence $\{ y_{i}\colon i \in \N\}$ of self-adjoint elements in $Y$ and $\epsilon >0$ choose a sequence $\{ x_{i} \colon i \in \N\} \subseteq X_{sa}$ such that $f(x_{i}) = y_{i}$ for each $i \in \N$. Since $X$ has weakly infinite bounded rank with respect to $K$, there exists a $K$-unessential sequence $\{ w_{i} \colon i \in \N\} \subseteq X$ such that $\| x_{i}-w_{i}\| < \epsilon$ for each $i \in \N$. We claim that $\{ z_{i} = f(w_{i}) \colon i \in \N\}$ is a $K$-unessential sequence in $Y$. Indeed, let $\delta > 0$ be a rational number. Since $\{ w_{i} \colon i \in \N\}$ is $K$-unessential in $X$, there exists a finite subset $D \subseteq \N$ and a $D$-tuple $\displaystyle\left( s_{i}\right)_{i \in D}$ such that
\begin{itemize}
\item[(i)$_{X}$]
$\| w_{i} - s_{i}\| \leq \delta$ for each $i \in D$;
\item[(ii)$_{X}$]
$\left\|\left(\sum_{i\in D}s_{i}^{2} \right)^{-1}\right\| \leq\displaystyle \frac{1}{K\cdot\delta^{2}}$.
\end{itemize}

Now consider the $D$-tuple $\left( r_{i} = f(s_{i})\right)_{ i \in D}$. Clearly
\begin{itemize}
\item[(i)$_{Y}$]
$\| z_{i} - r_{i}\| \leq \| f(w_{i})-f(s_{i})\| \leq \| f(w_{i}-s_{i})\|  \leq \| w_{i}-s_{i}\| \leq \delta$ for each $i \in D$;
\item[(ii)$_{Y}$]
$\left\|\left( \sum_{i\in D}r_{i}^{2}\right)^{-1}\right\| = \left\| f\left(\sum_{i\in D}s_{i}^{2} \right)^{-1}\right\|\leq \left\| \left(\sum_{i\in D}s_{i}^{2}\right)^{-1}\right\| \leq\displaystyle \frac{1}{K\cdot\delta^{2}}$.
\end{itemize}
Proof is completed.
\end{proof}

\begin{pro}\label{P:wcomparison}
Let $K > 0$. If the unital $C^{\ast}$-algebra $X$ has weakly infinite bounded rank with respect to $K$, then it has weakly infinite real rank.
\end{pro}
\begin{proof}
Let $\{ x_{i} \colon i \in \N\}$ be a sequence of self-adjoint elements in $X$ and let $\epsilon > 0$. Take a $K$-unessential sequence $\{ y_{i} \colon i \in \N\}$ in $X$ such that $\| x_{i} -y_{i}\| < \displaystyle\frac{\epsilon}{2}$ for each $i \in \N$. Since $\{ y_{i} \colon i \in \N\}$ is $K$-unessential, it contains a finite $K$-unessential subset $\{ y_{i} \colon i \in D\}$, $D \subseteq \N$. As in the proof of Proposition 4.6, there exists a $D$-tuple $( z_{i})_{i \in D}$ such that
\begin{itemize}
\item[(i)]
$\|y_{i}-z_{i}\| \leq \delta$ for each $i\in D$,
\item[(ii)]
$\sum_{i\in D}z_{i}^{2}$ is invertible,
\end{itemize}
Clearly, $\| x_{i}-z_{i}\| \leq \| x_{i}-y_{i}\| +\|y_{i}-z_{i}\| < \displaystyle\frac{\epsilon}{2}+\frac{\epsilon}{2} =\epsilon$, $i\in D$. According to (ii), $\sum_{i\in D}z_{i}^{2}$ is invertible, which shows that $X$ has weakly infinite real rank. 
\end{proof}

\begin{cor}\label{C:scomparison}
If a unital $C^{\ast}$-algebra has strongly infinite real rank,
then it has strongly infinite bounded rank with respect to any positive constant.
\end{cor}

\subsection{The commutative case}\label{SS:commutative}
If $X$ is a finite-dimensional compact space, then, according to Corollary \ref{C:rrcommutative}, $\operatorname{rr}(C(X))=\operatorname{br}_{1}(C(X)) =\dim X$ for any positive $K\leq 1$. Our next goal is to extend this result to the infinite-dimensional situation. 

First, recall that a compact Hausdorff space
$X$ is called {\em weakly infinite-\-dimen\-sional} \cite{ap:73} if for
any sequence 
$\{(F_i,H_i) \colon i \in \mathbb{N}\}$ of pairs of closed disjoint subsets of $X$ there
are partitions $L_i$ between $F_i$ and $H_i$ such that
$\bigcap_{i=1}^{\infty}L_i=\emptyset$. Here, $L_i\subset X$ is
called a partition between $F_i$ and $H_i$ if $L_i$ is closed   
in $X$ and $X\backslash L_i$ is decomposed as the union
$U_i\cup V_i$ of disjoint open sets with $F_i\subset U_i$
and $H_i\subset V_i$. Since $X$ is compact,
$\bigcap_{i=1}^{\infty}L_i=\emptyset$ is equivalent to
$\bigcap_{i=1}^{k}L_i=\emptyset$ for some $k\in\mathbb{N}$. If $X$ is not weakly infinite-dimensional, then it is {\em strongly infinite-dimensional}.

A standard example of a weakly infinite dimensional, but not finite-dimensional, metrizable  compactum can be obtained by taking the
one-point compactification $\alpha \left( \oplus\{ \mathbb{I}^{n} \colon n \in \N\}\right)$ of the discrete union of increasing dimensional cubes. The Hilbert cube $\mathbb{Q}$ is, of course, strongly infinite-dimensional.

\begin{thm}\label{T:iff}
Let $X$ be a compact Hausdorff space and $0<K\leq 1$. Then the following conditions are equivalent:
\begin{itemize}
\item[(a)]
$C(X)$ has weakly infinite bounded rank with respect to $K$;
\item[(b)]
$C(X)$ has weakly infinite real rank;
\item[(c)]
$X$ is weakly infinite-dimensional.
\end{itemize}
\end{thm}
\begin{proof}
(a)$\Longrightarrow$(b). This implication follows from Proposition \ref{P:wcomparison}
(which is valid for any -- not necessarily commutative -- unital $C^{\ast}$-algebras).

(b)$\Longrightarrow$(c).
Suppose that $C(X)$ has a weakly infinite real rank. Take an
arbitrary sequence $\{ (B_i,K_i) \colon i \in \mathbb{N}\}$ of pairs of disjoint
closed subsets of $X$ and define functions
$f_i\colon X\to [-1,1]$ such that $f_i(B_i)=-1$ and $f_i(K_i)=1$
for every $i \in \mathbb{N}$. Then, according to our hypothesis, there is
a sequence $\{g_i \colon i \in \mathbb{N}\}\subset C(X)$ of real-valued functions and an
integer $k$ with 
$\| f_i-g_i \| <1$, $i\in\mathbb{N}$, and
$\sum_{i=1}^{k}g_i^2(x)>0$ for each $x\in X$. If $C_i$ denotes
the set $g_i^{-1}(0)$, the last inequality means that
$\bigcap_{i=1}^{k}C_i=\emptyset$. Therefore, in order to prove
that $X$ is weakly infinite-dimensional, it only remains to
show each $C_i$ is a separator between $B_i$ and $K_i$. To
this end, we fix $i \in \mathbb{N}$ and observe that $\| f_i-g_i \| <1$ implies
the following inclusions: $g_i(B_i)\subseteq [-2,0)$, 
$g_i(K_i)\subseteq (0,2]$ and $g_i(X)\subseteq [-2,2]$. So,
$X\backslash C_i = U_i \cup V_i$, where $U_i=g_i^{-1}([-2,0))$ and
$V_i=g_i^{-1}((0,2])$. Moreover, $B_i\subseteq U_i$ and
$K_i\subseteq V_i$, i.e. $C_i$ separates $B_i$ and $K_i$.

(c)$\Longrightarrow$(a).
Let us show that the weak infinite-dimensionality of $X$
forces $C(X)$ to have a weakly infinite bounded rank with respect to $K$. To this end,
take any sequence $\{f_i \colon i \in \mathbb{N}\}\subset C(X)$ of real-valued functions and any positive number $\epsilon$.  It suffices to find another sequence $\{g_i \colon i \in \mathbb{N}\}$ of real-valued functions in
$C(X)$ such that $\|f_i-g_i\|\leq\epsilon$ for every $i \in\mathbb{N}$ and
$\sum_{i=1}^{m}g_i^2(x)>0$ for every
$x\in X$ and some $m\in\mathbb{N}$. Indeed,
if $\sum_{i=1}^{m}g_i^2(x)>0$ for every $x\in X$, then the function $\sum_{i=1}^{m}g_i^2$ is invertible. This, according to Proposition \ref{P:commutting}, is equivalent to the $K$-unessentiality of the $m$-tuple $(g_{1},\dots ,g_{m})$. 
On the other hand, $\sum_{i=1}^{m}g_i^2(x)>0$ for each $x\in X$ if and only if  $\bigcap_{i=1}^{m}g_i^{-1}(0)=\emptyset$. 
Further, since $X$ is compact,
the existence of  $m \in\mathbb{N}$ with $\bigcap_{i=1}^{m}g_i^{-1}(0)=\emptyset$ is equivalent to  
$\bigcap_{i=1}^{\infty}g_i^{-1}(0)=\emptyset$. Therefore, our
proof is reduced to constructing, for each $i \in\mathbb{N}$, a function $g_i$ 
which is $\epsilon$-close to $f_i$ and such that the
intersection of all $g_i^{-1}(0)$'s, $i\in\mathbb{N}$, is empty. 

For every $i \in \mathbb{N}$ let $c_i=\inf\{f_i(x):x\in X\}$ and
$d_i=\sup\{f_i(x):x\in X\}$. We can suppose, without loss
of generality, that each interval $(c_i,d_i)$ is not empty and contains $0$.  For every $i$ we choose $\eta_i>0$ such
that $\displaystyle\eta_i<\frac{\epsilon}{2}$ and  
$L_i=[-\eta_i, \eta_i ]\subset (c_i,d_i)$, $i\in\N$. Let $Q=\prod_{i=1}^{\infty}[c_i,d_i]$, 
$Q_0=\prod_{i=1}^{\infty}L_i$ be the topological products of all $[c_{i},d_{i}]$'s and $L_{i}$'s, respectively. Consider the diagonal product $f = \triangle \{ f_i \colon i\in\mathbb{N}\} \colon X \to Q$ and note that $H=\bigcap_{i=1}^{\infty}H_i$, where
$H = f^{-1}(Q_0)$ and $H_i=f_i^{-1}(L_i)$ for each $i \in \mathbb{N}$. We also consider the sets

 \[ F_i^{-}=f_i^{-1}([c_i,-\eta_i ]) \;\;
\text{and}\;\; F_i^{+}=f_i^{-1}([\eta_i , d_i]),\; i \in \mathbb{N} .\]

\noindent
Since $H$ is weakly infinite-dimensional (as a closed subset
of $X$), by \cite[Theorem 19, \S 10.4]{ap:73}, there
is a continuous map $p=(p_1,p_2,\dots )\colon H\to Q_0$ and a
pseudointerior point $b=\{ b_i \colon i \in \mathbb{N}\}\in Q_0$ (i.e. each $b_i$ lies
in the interior of the interval $L_i$)
such that

\[ b \not\in p(H),\;\; F_i^{-}\cap H\subset p_i^{-1}(\{ -\eta_i\} ),\;\;\text{and}\;\;
F_i^{+}\cap H\subset p_i^{-1}(\{\eta_i\} ),\; i\in\mathbb{N} .\]

\noindent
Since each $b_i$ is an interior point of
$L_i=[-\eta_i , \eta_i ]$, there exists homeomorphisms
$s_i\colon L_i\to L_i$ which leaves the endpoints $-\eta_i$ and $\eta_i$ fixed and such that $s_i(b_i)=0$.
Let $s= \triangle \{ s_i \colon i \in \mathbb{N}\} \colon Q_0\to Q_0$ and $q=s\circ p$. Obviously $s(b)=\bf{0}$ and ${\bf{0}}\not\in q(H)$, where $\bf{0}$ denotes the point of $Q_0$ having all coordinates $0$. Further observe that if $q_{i} = \pi_{i}\circ q$, where $q_{i} \colon Q_{0} \to L_{i}$ denotes the natural projection onto the $i$-th coordinate, then

\[ F_i^{-}\cap H\subset q_i^{-1}(\{ -\eta_i \})\;\;\text{and}\;\;
F_i^{+}\cap H\subset q_i^{-1}(\{ \eta_i \}) ,\;\; i \in \mathbb{N}.\]

\noindent
Therefore, each $q_i$, $i \in \mathbb{N}$, is a function from $H$ into $L_i$
satisfying the following condition:
 $q_i(F_i^{-}\cap H)=f_i(F_i^{-}\cap H_i)=-\eta_i$
and  $q_i(F_i^{+}\cap H)=f_i(F_i^{+}\cap H_i)=\eta_i$.
Let $h_i\colon H_i\to L_i$ be an extension of $q_{i}$, $i \in \mathbb{N}$.
Note that the restrictions of $h_{i}$ and $f_{i}$ onto the sets $F_i^{-}\cap H_i$ and $F_i^{+}\cap H_i$ coincide.
Finally, define $g_i\colon X\to [c_i,d_i]$ by letting

\[ g_{i}(x) = 
\begin{cases}
h_{i}(x) ,\;\;\text{if}\;\; x \in H_{i};\\
f_{i}(x) ,\;\;\text{if}\;\; x \in X-H_{i}.\\
\end{cases}
\]

To finish the proof of the "if" part, we need to show that 
$g_i(x)$ is $\epsilon$-close to $f_i(x)$ for each $i \in \mathbb{N}$ and
$x\in X$, and that $\bigcap_{i=1}^{\infty}g_i^{-1}(0)=\emptyset$.
Since $g_i$ and $f_i$ are identical outside $H_i$, the first
condition is satisfied for $x\not\in H_i$. If $x\in H_i$,
then both $f_i(x)$ and $g_i(x)$ belong to $L_i$, so again
$|f_i(x)-g_i(x)| < \epsilon$. To prove the second condition,
observe first that $x\not\in H$ implies $x\not\in H_j$ for
some $j$. Hence, $g_j(x)=f_j(x)\not\in L_j$, so $g_j(x)\neq 0$.
If $x\in H$, then $g_i(x)=q_i(x)$ for all $i$ and, because
${\bf{0}}\not\in q(H)$, at least one $g_i(x)$ must be different
from $0$. Thus, $\bigcap_{i=1}^{\infty}g_i^{-1}(0)=\emptyset$. 
\end{proof}

Let $C^{\ast}\left( {\mathbb F}_{\infty}\right)$ denote the group
$C^{\ast}$-algebra  
of the free group on countable number of generators. 
It is clear that
$\operatorname{rr}\left( C^{\ast}\left(
{\mathbb F}_{\infty}\right)\right) > n$ for each $n$. Our results
imply much stronger observation.

\begin{cor}\label{C:notweakly}
The group
$C^{\ast}$-algebra $C^{\ast}\left( {\mathbb F}_{\infty}\right)$ 
of the free group on countable number of generators
has strongly infinite real rank.
\end{cor}
\begin{proof}
It is well known that every separable unital $C^{\ast}$-algebra
is an image of
$C^{\ast}\left( {\mathbb F}_{\infty}\right)$
under a surjective ${\ast}$-homomorphism. In particular, there
exists a surjective ${\ast}$-homo\-mor\-phism
$f \colon C^{\ast}\left( {\mathbb F}_{\infty}\right) \to C(Q)$,
where $Q$ denotes the Hilbert cube. It is well known (see,
for instance, \cite[\S 10.5]{ap:73}) that the Hilbert cube $Q$ 
is strongly infinite dimensional.
By Theorem \ref{T:iff}, $C(Q)$ has strongly infinite real rank. Finally, by Proposition \ref{L:homomorphism},
real rank of $C^{\ast}\left( {\mathbb F}_{\infty}\right)$ must also be strongly infinite.
\end{proof}

\begin{pro}\label{P:two}
Let $X$ and $Y$ be unital $C^{\ast}$-algebras with weakly infinite real rank. Then $X \oplus Y$ also has weakly infinite real tank.
\end{pro}
\begin{proof}
Let $\{ (x_{i}, y_{i}) \colon i \in \mathbb{N}\}$ be a sequence of self-adjoint elements of $X\oplus Y$ and $\epsilon > 0$. Since both $X$ and $Y$ have weakly infinite reak rank there exist sequences $\{ z_{i} \colon i \in \mathbb{N}\}$ and $\{ w_{i} \colon i \in \mathbb{N}\}$ of self-adjoint elements of $X$ and $Y$ respectively such that
\begin{itemize}
\item[(i)]
$\| x_{i} -z_{i}\| < \epsilon$ for every $i \in \mathbb{N}$;
\item[(ii)]
$\| y_{i} -w_{i}\| < \epsilon$ for every $i \in \mathbb{N}$;
\item[(iii)]
for some $n \geq 1$ the element $\sum_{i=1}^{n}z_{i}^{2}$ is invertible;
\item[(iv)]
for some $m \geq 1$ the element $\sum_{i=1}^{m}w_{i}^{2}$ is invertible.
\end{itemize}

Without loss of generality we may assume that $ n \geq m$. Obviously, according to (iv), $\sum_{i=1}^{n}w_{i}^{2}$ is also invertible.
Next consider the sequence $\{ (z_{i},w_{i}) \colon i \in \mathbb{N}\}$.
Note that $\| (x_{i},y_{i}) - (z_{i},w_{i})\| = \| (x_{i}-z_{i},y_{i}-w_{i})\| = \max\{ \| x_{i}-z_{i}\| ,\| y_{i}-w_{i}\|\} <\epsilon$ and that, according to (iii) and (iv),
$\sum_{i=1}^{n}(z_{i},w_{i})^{2} = \sum_{i=1}^{n}(z_{i}^{2},w_{i}^{2})$ is also invertible. 
\end{proof}

Next statement provides a formal example of a unital $C^{\ast}$-algebra
of weakly infinite, but not finite real rank. 
\begin{pro}\label{P:example}
Let $X = \alpha \left(\oplus\{ I^{n}\colon n \in \mathbb{N}\}\right)$ be the one-point compactification of the discrete topological sum of increasing-dimensional cubes. In other words, $C(X) = \prod\{ C(I^{n})\colon n \in \mathbb{N}\}$ (here $\prod$ stands for the direct product of indicated $C^{\ast}$-algebras).
Then $C(X)$ has weakly infinite, but not finite real rank.
\end{pro}
\begin{proof}
Obviously $X$ is countably dimensional and hence, by \cite[Corollary 1, \S 10.5]{ap:73}, it is weakly infinite dimensional. By Theorem \ref{T:iff}, $C(X)$ has weakly infinite real rank. It only remains to note that $\operatorname{rr}(X) > n$ for any $n \in \mathbb{N}$.
\end{proof}

In conclusion let us note that there exist non-commutative $C^{\ast}$-algebras with similar properties (compare with Corollary \ref{C:notweakly}).
\begin{cor}\label{C:nonabelian}
There exist non-commutative untal $C^{\ast}$-algebras of weakly infinite, but not finite real rank.
\end{cor}
\begin{proof}
Let $X$ be as in Proposition \ref{P:example} and $A$ be a non-commutative unital $C^{\ast}$-algebra of a finite real rank. According to Propositions \ref{P:finite} and \ref{P:two}, the product $C(X) \oplus A$ has weakly infinite real rank, It is clear that $C(X) \oplus A$ is non-commutative and does not have a finite real rank.
\end{proof}


\end{document}